\documentclass[A4,12pt]{amsart}
\usepackage{amsfonts,amsmath,amssymb,amsthm}
\usepackage{xcolor}

\numberwithin{equation}{section}
\newtheorem{theorem}{Theorem}[section]
\newtheorem{lemma}{Lemma}[section]
\newtheorem{proposition}{Proposition}[section]
\newtheorem{corollary}{Corollary}[section]

\newtheorem{remark}{Remark}[section]

\newtheorem{definition}{Definition}[section]
\newtheorem{example}{Example}[section]
\numberwithin{equation}{section}
\numberwithin{equation}{section}
\begin{document}
\title{Isotropic Riemannian Maps and Helices along Riemannian Maps}
\date{}
\author{G\"{o}zde \"{O}zkan T\"{u}kel}
\address[G. \"{O}. T\"{u}kel]{Isparta University of Applied Sciences, Faculty of Technology, Department of Engineering Basic Sciences,
Isparta, Turkey}
\email[G. \"{O}. T\"{u}kel]{gozdetukel@isparta.edu.tr}
\author{Bayram \c{S}ahin}
\address[B. \c{S}ahin]{Ege University, Faculty of Science, Department of
Mathematics, \.{I}zm\i r-Turkey}
\email[B. \c{S}ahin]{bayram.sahin@ege.edu.tr}
\author{Tunahan Turhan}
\address[T. Turhan]{S\"{u}leyman Demirel University, Faculty of Education, Division of Elementary Mathematics, Isparta,
Turkey}
\email[T. Turhan]{tunahanturhan@sdu.edu.tr}
\maketitle

$\frac{{}}{{}}$\noindent \textbf{Abstract}. This work has two main purposes.  The first aim is to study isotropic Riemannian maps as a generalization of isotropic immersions.  For this purpose, the concept of isotropic Riemannian map is presented, an example is given and a characterization is obtained.  The second aim is to study the helices along Riemannian map. For this purpose, by using the notion of isotropic Riemannian map and the notion of helix on the manifold, a characterization is obtained for the transportation of a helix on the total manifold to the target manifold along a Riemannian map.


\newline

\section{\textbf{Introduction}}
In differential geometry, introducing and using appropriate transformations to compare two geometric objects is one of the main features. In the theory of manifolds, immersions are the most used maps in this sense.
For example, Nomizu and Yano in \cite{Nomizu} showed that when a circle on
the submanifolds is carried along the immersion to the ambient manifold, such
submanifolds are umbilical and their mean curvature vector field is
parallel. S. Maeda in \cite{Maeda} has shown that if a circle on the
submanifolds is carried to curve with the constant first curvature on the
manifold by immersion, then the submanifold is the isotropic submanifold
introduced in \cite{O'Neil}. Also, Ikawa in \cite{Ikawa} generalized Maeda's results in  \cite{Maeda} by using helices. Later this result had been extended to the semi-Riemannian case \cite{Ikawa2} by same author, see also:, \cite{Ekmekci1}, \cite{Ekmekci2}, \cite{EH}, \cite{Ekmekci-Ilarslan}, \cite{Ilarslan}, \cite{Song}. Geometric properties and characterizations of
isometric immersions have been investigated by many researches, see: \cite{Chen}.

The basic properties of Riemannian submersions were firstly given by Gray \cite{Gray} and O'Neill \cite{O'Neill}. Riemannian submersions and isometric immersions are special maps between Riemannian manifolds. So, Fischer \cite{Fischer} defined Riemannian maps which is a generalization of them in $1992$ as follows. Assume that $\mathcal{T}:(\mathfrak{M}_1,g_{_{\mathfrak{M}_1}})\rightarrow
(\mathfrak{M}_2,g_{_{\mathfrak{M}_2}})$ is a $ C^{\infty} $ map from the Riemannian manifold $\mathfrak{M}_1 $ with $\dim \mathfrak{M}_1=m$ to the Riemannian manifold $\mathfrak{M}_2 $ with $\dim \mathfrak{M}_2=n$, where $%
0<rank\mathcal{T}<\min \{m,n\}.$ Thus, we represent
the kernel space of $\mathcal{T}_{\ast }$ by $\ker \mathcal{T}_{\ast }$ and $(\ker \mathcal{T}_{\ast })^{\perp }$ is orthogonal complementary space to $\ker
\mathcal{T}_{\ast }.$ So, we have%
\begin{equation*}
T\mathfrak{M}_1=\ker \mathcal{T}_{\ast }\oplus (\ker \mathcal{T}_{\ast })^{\perp },
\end{equation*}%
where $ T\mathfrak{M}_1 $ is the tangent bundle of $\mathfrak{M}_1$. $range\mathcal{T}_{\ast }$ denotes the range of $\mathcal{T}_{\ast }$  and  $(range\mathcal{T}_{\ast })^{\perp }$ denotes the
orthogonal complementary space to $range\mathcal{T}_{\ast
} $ in $T\mathfrak{M}_2$. The tangent bundle $T\mathfrak{M}_2$ of $%
\mathfrak{M}_2$ is given by%
\begin{equation*}
T\mathfrak{M}_2=(range\mathcal{T}_{\ast })\oplus (range\mathcal{T}_{\ast })^{\perp }.
\end{equation*}%
Now, a smooth map $\mathcal{T}:(\mathfrak{M}_1^{m},g_{_{\mathfrak{M}_1}})\rightarrow (\mathfrak{M}_2^{n},g_{_{\mathfrak{M}_2}})$ is called
Riemannian map if $\mathcal{T}_{\ast }$ satisfies $g_{_{\mathfrak{M}_2}}(\mathcal{T}_{\ast
}X_{1},\mathcal{T}_{\ast }X_{2})=g_{_{M}}(X_{1},X_{2}),$ for $X_{1}$, $X_{2}$ vector fields tangent to H \cite%
{Fischer}. 

Comparing immersions with Riemann submersions and Riemann maps, there are few results for the characterization of such transformations with the help of curves.  Indeed, there are no other studies in the literature, except for the results regarding the geodesics at reference \cite{FIP} for Riemann submersions and the results about circles at reference \cite{Sahin4} for Riemannian maps. Since curves are the basic notion of differential geometry and their useful role in immersion theory is taken into account, obtaining new results on the characterization of Riemannian submersions and Riemannian maps by means of curves is a subject that needs to be investigated.

Also, isotropic  immersions were defined by O'Neill \cite{O'Neil} and later it has been shown that this notion is an important tool in geometric characterization \cite{Maeda}.  The concept of isotropic submersions has been introduced by \c{S}ahin and Erdogan in \cite{Erdogan-Sahin}. As far as we know, the notion of  isotropic Riemannian maps has not yet been introduced in the literature.

We first define the concept of an isotropic Riemannian map
and  obtain a condition for isotropicity of a
Riemannian map in terms of second fundamental form. Then by using
characterization of an ordinary helix in a Riemannian manifold, we
generalize Ikawa's theorem.

In the second section, the basic notions to be used in the paper are presented. First, the concept of helix on the manifold is given. Then, Gauss and Weingarten formulas are introduced for Riemannian maps with connections defined along a map. These notions defined along a map are presented in some detail, especially considering that the readers who study the submanifold theory may not be familiar with these notions. In the third section, the concept of isotropic Riemannian map is presented and characterization of such maps is obtained. We give the following result in the last section; if there is a helix on the source manifold of a Riemannian map, we investigate what this property tells  for the Riemannian map when a helix is transformed by the Riemannian map on the target manifold.

\section{\textbf{Preliminaries}}
A regular curve $\alpha =\alpha (s)$ parametrized by arc lengh $s$ is called
an ordinary helix if there exist unit vector fields $V_2$ and $V_3$ along $\alpha
$\ such that%
\begin{equation*}
\nabla _{V_1}V_1=\kappa V_2,\text{ \ }\nabla _{V_1}V_2=-\kappa V_1+\tau V_3,\text{ \ }%
\nabla _{V_1}V_1=-\tau V_2,
\end{equation*}%
where $V_1$ denotes the tangent vector field of $\alpha $, $\kappa $ is the curvature and $\tau
$ is the torsion of $\alpha$. An ordinary helix satisfies the following equation
\begin{equation}
\nabla^{3} _{V_1}V_1+K^2\nabla _{V_1}V_1=0,
\label{helix}
\end{equation}%
where $ K^2 $ is a constant. If $ \tau=0 $, then helix reduces to the circle \cite{Ikawa}.

Assume that $(\mathfrak{M}_1,g_{_{\mathfrak{M}_1}})$ and $(\mathfrak{M}_2,g_{_{\mathfrak{M}_2}})$ are Riemannian manifolds, $\mathcal{T}:(\mathfrak{M}_1,g_{_{\mathfrak{M}_1}})\rightarrow (\mathfrak{M}_2,g_{_{\mathfrak{M}_2}})$ a smooth map between them and $ \gamma $ a curve on $\mathfrak{M}_1 $. $ \gamma $ is called a horizontal curve if $ \dot{\gamma}(t)\in (ker \mathcal{T}_{*})^\bot$ for any $ t\in I $. If $ \gamma $ is a helix with $ \dot{\gamma}(t)\in (ker \mathcal{T}_{*})^\bot$ for any $ t\in I $, then it is called as horizontal helix.

Let $\mathcal{T}$ be a Riemannian map between the manifolds  $(\mathfrak{M}_1,g_{_{\mathfrak{M}_1}})$ and $(\mathfrak{M}_2,g_{_{\mathfrak{M}_2}})$, $p_{_{2}}=\mathcal{T}(p_{1})$ for each $p_{1}\in \mathfrak{M}_1$. Suppose that $\nabla ^{\mathfrak{M}_2}$ and $\nabla ^{\mathfrak{M}_1}$ represent
the connections on $(\mathfrak{M}_2,g_{_{\mathfrak{M}_2}})$ and $(\mathfrak{M}_1,g_{_{\mathfrak{M}_1}}),$ respectively. The second fundamental
form of $\mathcal{T}$ can be given as follows%
\begin{equation}
(\nabla \mathcal{T}_{\ast })(X_{1},X_{2})=\overset{\mathfrak{M}_2}{\nabla ^{\mathcal{T}}}_{X_{1}}\mathcal{T}_{\ast }(X_{2})-\mathcal{T}_{\ast
}(\nabla _{X_{1}}^{\mathfrak{M}_1}X_{2})
\label{2.1}
\end{equation}%
for $X_{1}$, $X_{2}\in \Gamma (T\mathfrak{M}_1),$ where $\overset{\mathfrak{M}_2}{\nabla ^{\mathcal{T}}}$ is the
pullback connection of $\nabla ^{\mathfrak{M}_2}.$ For $\forall X_{1},X_{2}\in \Gamma ((\ker \mathcal{T}_{\ast
p_{1}})^{\perp }),$ $ (\nabla \mathcal{T}_{\ast }) $ is symmetric and has no components in $range\mathcal{T}_{\ast
}. $ So, we can write the following
\begin{equation}
g_{_{N}}((\nabla \mathcal{T}_{\ast })(X_{1},X_{2}),\mathcal{T}_{\ast }(X_{3}))=0,
\label{2.2}
\end{equation}
$\forall X_{1},X_{2},X_{3}\in
\Gamma ((\ker \mathcal{T}_{\ast p_{1}})^{\perp })$ \cite{Sahink}.

Now we give some basic formulas for Riemannian maps defined from the total manifold $(\mathfrak{M}_1,g_{_{\mathfrak{M}_1}})$ to  the target manifold $%
(\mathfrak{M}_2,g_{_{\mathfrak{M}_2}}) $. For $%
X_{1},X_{2}\in \Gamma ((\ker \mathcal{T}_{\ast p_{1}})^{\perp })$ and $U_{1}\in \Gamma
((range\mathcal{T}_{\ast })^{\perp }),$ we have%
\begin{equation}
\overset{\mathfrak{M}_2}{\nabla^{\mathcal{T}}}_{X_{1}}U_{1}=-S_{U_{1}}\mathcal{T}_{\ast }X_{1}+\nabla
_{X_{1}}^{\mathcal{T}^{\perp }}U_{1},  \label{2.3}
\end{equation}%
where $S_{X_{1}}\mathcal{T}_{\ast }X_{1}$ is the tangential component of $\overset{\mathfrak{M}_2}{\nabla^{\mathcal{T}}}_{X_{1}}U_{1}$ and we have%
\begin{equation}
g_{_{\mathfrak{M}_2}}(S_{U_{1}}\mathcal{T}_{\ast }X_{1},\mathcal{T}_{\ast }X_{2})=g_{_{\mathfrak{M}_2}}(U_{1},(\nabla \mathcal{T}_{\ast })(X_{1},X_{2})).
\label{2.4}
\end{equation}%
Since $(\nabla \mathcal{T}_{\ast })$ is symmetric, $S_{U_{1}}$ is a symmetric linear transformation of $range\mathcal{T}_{\ast }.$ On the other hand, we have the following covariant derivatives%
\begin{equation}
\left( \tilde{\nabla}_{X_{1}}(\nabla \mathcal{T}_{\ast })\right) (X_{2},X_{3})=\nabla
_{X_{1}}^{\mathcal{T}^{\perp }}(\nabla \mathcal{T}_{\ast })(X_{2},X_{3})-(\nabla \mathcal{T}_{\ast })(\overset{\mathfrak{M}_1}{%
\nabla }_{X_{1}}X_{2},X_{3})-(\nabla \mathcal{T}_{\ast })(X_{2},\overset{\mathfrak{M}_1}{\nabla }_{X_{1}}X_{3})
\label{2.5}
\end{equation}%
and%
\begin{equation}
\left( \tilde{\nabla}_{X_{1}}S\right) _{U_{1}}\mathcal{T}_{\ast }(X_{2})=\mathcal{T}_{\ast }(\overset{\mathfrak{M}_1}{%
\nabla }_{X_{1}}\text{~}^{\ast }\mathcal{T}_{\ast }(S_{U_{1}}\mathcal{T}_{\ast }(X_{2})))-S_{_{\left( \nabla
_{X_{1}}^{\mathcal{T}^{\perp }}\right) U_{1}}}\mathcal{T}_{\ast }(X_{2})-S_{U_{1}}\overset{\mathfrak{M}_2~~~}{P\nabla _{X_{1}}^{\mathcal{T}}}\mathcal{T}_{\ast
}(X_{2}),  \label{2.6}
\end{equation}%
where $^{\ast }\mathcal{T}_{\ast }$ denotes the projection morphism on $range\mathcal{T}_{\ast }$
and $^{\ast }\mathcal{T}_{\ast }$ is the adjoint map of $\mathcal{T}_{\ast }$ \cite{Sahink}. In
the following lemma we give a relation obtained from (\ref{2.5}) and (\ref{2.6}).

\begin{lemma}
Let $(\mathfrak{M}_1,g_{_{\mathfrak{M}_1}})$,
$(\mathfrak{M}_2,g_{_{\mathfrak{M}_2}})$ be Riemannian manifolds and $\mathcal{T}$  a Riemannian map between them. For  $\forall X_{1},X_{2}\in \Gamma ((\ker \mathcal{T}_{\ast
p_{1}})^{\perp })$ and $U_{1}\in \Gamma ((range\mathcal{T}_{\ast })^{\perp }),$ we have%
\begin{equation}
<\left( \tilde{\nabla}_{X_{1}}(\nabla \mathcal{T}_{\ast })\right) (X_{2},X_{3}),U_{1}>=<\left( \tilde{%
\nabla}_{X_{1}}S\right) _{U_{1}}\mathcal{T}_{\ast }(U_{1}),\mathcal{T}_{\ast }(X_{3})>.  \label{2.7}
\end{equation}
\end{lemma}

\begin{proof}
Taking inner product (\ref{2.5}) with $U_{1},$ we have%
\begin{equation*}
\begin{array}{ccl}
<\left( \tilde{\nabla}_{X_{1}}(\nabla \mathcal{T}_{\ast })\right) (X_{2},X_{3}),U_{1}> & = & <\nabla
_{X_{1}}^{\mathcal{T}^{\perp }}(\nabla \mathcal{T}_{\ast })(X_{2},X_{3}),U_{1}>-<(\nabla \mathcal{T}_{\ast })(\overset{\mathfrak{M}_1}{\nabla }_{X_{1}}X_{2},X_{3}),U_{1}> \\
&  & -<(\nabla \mathcal{T}_{\ast })(X_{2},\overset{\mathfrak{M}_1}{\nabla }_{X_{1}}X_{3}),U_{1}>.%
\end{array}%
\end{equation*}%
If we take inner product $(2.7)$ with $\mathcal{T}_{\ast }(X_{3}),$ we obtain%
\begin{equation*}
\begin{array}{ccl}
<\left( \tilde{\nabla}_{X_{1}}S\right) _{U_{1}}\mathcal{T}_{\ast }(X_{2}),\mathcal{T}_{\ast }(X_{3})> & = &
<\mathcal{T}_{\ast }(\overset{\mathfrak{M}_1}{\nabla }_{X_{1}}\text{~}^{\ast }\mathcal{T}_{\ast }(S_{U_{1}}\mathcal{T}_{\ast
}(X_{2}))),\mathcal{T}_{\ast }(X_{3})> \\
&  & -<S_{_{\left( \nabla
_{X_{1}}^{\mathcal{T}^{\perp }}\right) U_{1}}}\mathcal{T}_{\ast
}(X_{2}),\mathcal{T}_{\ast }(X_{3})>-<S_{U_{1}}\overset{\mathfrak{M}_2}{P\nabla _{X_{1}}^{\mathcal{T}}}\mathcal{T}_{\ast
}(X_{2}),\mathcal{T}_{\ast }(X_{3})>.%
\end{array}%
\end{equation*}%
Using (\ref{2.4}), we can write the  following equalites%
\begin{equation}
<(\nabla \mathcal{T}_{\ast })(\overset{\mathfrak{M}_1}{\nabla }_{X_{1}}X_{2},X_{3}),U_{1}>=<S_{U_{1}}\overset{\mathfrak{M}_2}{P\nabla _{X_{1}}^{\mathcal{T}}}\mathcal{T}_{\ast }(X_{2}),X_{3}>  \label{2.8}
\end{equation}%
and%
\begin{equation*}
<S_{U_{1}}\mathcal{T}_{\ast }(X_{2}),\mathcal{T}_{\ast }(X_{3})>=<(\nabla \mathcal{T}_{\ast })(X_{2},X_{3}),U_{1}>.
\end{equation*}%
If we take derivative of the last equation and use (\ref{2.2}), we get%
\begin{equation}
\begin{array}{c}
<\mathcal{T}_{\ast }(\overset{\mathfrak{M}_1}{\nabla }_{X_{1}}\text{~}^{\ast }\mathcal{T}_{\ast }(S_{U_{1}}\mathcal{T}_{\ast
}(X_{2}))),\mathcal{T}_{\ast }(X_{3})>+<(\nabla \mathcal{T}_{\ast })(X_{2},\overset{\mathfrak{M}_1}{\nabla }_{X_{1}}X_{3}),U_{1}> \\
=<\nabla _{X_{1}}^{\mathcal{T}^{\perp }}(\nabla \mathcal{T}_{\ast })(X_{2},X_{3}),U_{1}>+<S_{_{\left( \nabla
_{X_{1}}^{\mathcal{T}^{\perp }}\right) U_{1}}}\mathcal{T}_{\ast
}(X_{2}),\mathcal{T}_{\ast }(X_{3})>.%
\end{array}
\label{2.9}
\end{equation}

If we take into consideration (\ref{2.8}) and (\ref{2.9}), we have (\ref{2.7}).
\end{proof}

We consider that $\mathcal{T}$ is a Riemannian map from a connected Riemannian manifold $%
(\mathfrak{M}_1,g_{_{\mathfrak{M}_1}}),$ $dim \mathfrak{M}_1 \geq 2$ to a Riemannian manifold $(\mathfrak{M}_2,g_{_{\mathfrak{M}_2}}).$ We know
that $\mathcal{T}$ is an umbilical Riemannian map at $p_{1}\in \mathfrak{M}_1,$ if the following is satisfied%
\begin{equation*}
S_{U_{1}}\mathcal{T}_{\ast p_{1}}((X_{1})_{p_{1}})=g_{_{\mathfrak{M}_2}}(H_{2},U_{1})\mathcal{T}_{\ast p_{1}}((X_{1})_{p_{1}})
\end{equation*}%
for $X_{1}\in \Gamma (range\mathcal{T}_{\ast })$, $U_{1}\in \Gamma ((range\mathcal{T}_{\ast })^{\perp })$
and$\ H_{2}\in (range\mathcal{T}_{\ast })^{\perp }$. If $\mathcal{T}$ is umbilical for $\forall
p_{1}\in \mathfrak{M}_1$, we know that $\mathcal{T}$ is umbilical \cite{Sahink}%
.

\section{Isotropic Riemannian maps}
Now, we present the notion of isotropic Riemannian map and get a characterization for such maps.
\begin{definition}
\textbf{\ (h-Isotropic Riemannian map).} A Riemannian map $\mathcal{T}:\mathfrak{M}_1 \rightarrow \mathfrak{M}_2$
is said to be h-isotropic at $p\in \mathfrak{M}_1 $ if $\lambda (X_{1})=\left\Vert (\nabla
\mathcal{T}_{\ast })(X_{1},X_{1})\right\Vert /\left\Vert \mathcal{T}_{\ast }X_{1}\right\Vert ^{2}$ doesn't
depend upon the selection of $X_{1}\in \Gamma ((\ker \mathcal{T}_{\ast })^{\perp }.$ If the map
is h-isotropic at all points, the map is called as h-isotropic. Also, if $\lambda =\lambda (p)$ is constant along $\mathcal{T}$,  $\mathcal{T}$ is called a constant $(\lambda h-)$isotropic map.

\begin{remark}
It is easy to see that the above notion is generalization of isotropic immersion, however there is no inclusion relation between isotropic Riemannian submersions and $(h-)$isotropic Riemannian maps.
\end{remark}

\begin{proposition}
Let $ \varphi : \mathfrak{M}_1^{m} \rightarrow  \mathfrak{M}_2^{n}  $ be a Riemannian submersion and $  \psi  : \mathfrak{M}_2^{n} \rightarrow  \mathfrak{M}_3^{k} $ a $ \lambda- $isotropic immersion. Then the Riemannian map $ \mathcal{T}=\psi \circ \varphi $ is a $ \lambda h- $ isotropic map.
\end{proposition}

\begin{proof}
We consider that $ \mathcal{T}=\psi \circ \varphi $ is a Riemannian map, where $ \varphi : \mathfrak{M}_1^{m} \rightarrow  \mathfrak{M}_2^{n}  $ is a Riemannian submersion and $  \psi  : \mathfrak{M}_2^{n} \rightarrow  \mathfrak{M}_3^{k} $ is a $ \lambda- $isotropic immersion. Then we have for $ X_{1},X_{2}\in (ker\mathcal{T}_{*})^\perp $
\begin{equation}
(\nabla (\psi \circ \varphi )_{*})(X_{1},X_{2})=\psi_{*}((\nabla \varphi_{*})(X_{1},X_{2}))+(\nabla \psi_{*})(\varphi_{*}(X_{1}),\varphi_{*}(X_{2})).
 \label{aa}
\end{equation}
On the other hand, we get
\begin{equation}
g_{\mathfrak{M}_1}((\psi \circ \varphi)_{*} (X_{1}),(\psi \circ \varphi)_{*}(X_{2}))=g_{\mathfrak{M}_2}( \psi_{*} (X_{1}), \psi_{*}(X_{2}))=g_{\mathfrak{M}_3}( X_{1}, X_{2}).
\label{bb}
\end{equation}
From ( \ref{aa} ), we obtain
\begin{equation}
\begin{array}{cc}
\Vert (\nabla (\psi \circ \varphi )_{*})(X_{1},X_{1}) \Vert ^{2} =  \Vert \psi_{*}((\nabla \varphi_{*})(X_{1},X_{1})) \Vert ^{2} \\ +2g_{\mathfrak{M}_3}(\psi_{*}((\nabla \varphi_{*})(X_{1},X_{1})),(\nabla \psi_{*})(\varphi_{*}(X_{1}),\varphi_{*}(X_{1}) )  \\  + \Vert (\nabla \psi_{*})(\varphi_{*}(X_{1}),\varphi_{*}(X_{1}) \Vert ^2.
\end{array}
\label{cc}
\end{equation}
Because $  \psi $ is a $ \lambda- $isotropic immersion and due to Riemannian submersion $ \varphi $, $(\nabla \varphi_{*})(X_{1},X_{1})=0$ for $ X_{1}\in (ker\mathcal{T}_{*})^\perp $, (\ref{cc}) and (\ref{bb}) show that $ \mathcal{T} $ is a $\lambda h-$ isotropic map.
\end{proof}
\emph{The following example can be given as an application of above result.}
\begin{example}
Let $\mathcal{T}$ be a Riemannian map given by
\begin{equation*}
\begin{array}{ccccc}
\mathcal{T} & : & \mathbb{R}^{4} & \rightarrow & \mathbb{R}^{3} \\
&  & (x_{1},x_{2},x_{3},x_{4}) & \rightarrow & (\dfrac{(x_{1}-x_{2})^2}{2}-x_{3}^2%
,\sqrt{2}(x_{1}-x_{2})x_{3},0).%
\end{array}%
\end{equation*}%
Since
\begin{equation*}
\begin{array}{ccccc}
\varphi & : & \mathbb{R}^{4} & \rightarrow & \mathbb{R}^{2}
\\
&  & (x_{1},x_{2},x_{3},x_{4}) & \rightarrow & (\frac{x_{1}-x_{2}}{\sqrt{2}},x_{3})
\end{array}
\end{equation*}
is a Riemannian submersion and
\begin{equation*}
\begin{array}{ccccc}
\psi & : & \mathbb{R}^{2} & \rightarrow & \mathbb{R}^{3}
\\
&  & (x_{1},x_{2}) & \rightarrow & (x_{1}^{2}-x_{2}^{2},2x_{1}x_{2},0),
\end{array}
\end{equation*}
where $(x_{1}^{2}+x_{2}^{2})^2=1  $ (see, \cite{O'Neilll}), a $(\lambda -)$ isotropic immersion, the Riemannian map $ \mathcal{T}=\psi \circ \varphi $ is a  $(\lambda h-)$%
isotropic Riemannian map.
\end{example}
\end{definition}
The following lemma gives a criteria for a  h-isotropicity Riemannian map.
\begin{lemma} \label{Lemma}
Let $\mathcal{T}:\mathfrak{M}_1\rightarrow \mathfrak{M}_2$ be a Riemannian map. $\mathcal{T}$ is h-isotropic at $p\in \mathfrak{M}_1$
iff the second fundamental form $\nabla \mathcal{T}_{\ast }$ satisfies%
\begin{equation}
g_{_{\mathfrak{M}_2}}((\nabla \mathcal{T}_{\ast })(X_{1},X_{1}),(\nabla \mathcal{T}_{\ast })(X_{1},X_{2}))=0  \label{3.2}
\end{equation}%
for an arbitrary orthogonal couple $X_{1},X_{2}\in \Gamma ((\ker \mathcal{T}_{\ast
p_{1}})^{\perp }).$
\end{lemma}
\begin{proof}
Let $f:\Gamma ((\ker \mathcal{T}_{\ast p_{1}})^{\perp })\rightarrow \mathbb{R}$ be a
quadrilinear function such that for $\forall x_{1}$, $x_{2}$, $u_{1}$, $u_{2}\in \Gamma
((\ker \mathcal{T}_{\ast p_{1}})^{\perp })$
\begin{equation*}
f(x_{1},x_{2},u_{1},u_{2})=g_{_{\mathfrak{M}_2}}((\nabla \mathcal{T}_{\ast })(x_{1},x_{2}),(\nabla \mathcal{T}_{\ast })(u_{1},u_{2}))-\lambda
^{2}g_{_{\mathfrak{M}_2}}(\mathcal{T}_{\ast }(x_{1}),\mathcal{T}_{\ast }(x_{2}))g_{_{N}}(\mathcal{T}_{\ast }(u_{1}),\mathcal{T}_{\ast }(u_{2})).
\end{equation*}%
If a Riemannian map $\mathcal{T}$ is $\lambda h-$isotropic, we have for $\forall u_{3}\in
\Gamma ((\ker \mathcal{T}_{\ast p_{1}})^{\perp })$%
\begin{equation*}
\begin{array}{ccl}
B(u_{3}) & = & f(u_{3},u_{3},u_{3},u_{3}), \\
& = & g_{_{N}}((\nabla \mathcal{T}_{\ast })(u_{3},u_{3}),(\nabla \mathcal{T}_{\ast })(u_{3},u_{3}))-\lambda
^{2}g_{_{M}}(u_{3},u_{3})g_{_{M}}(u_{3},u_{31}), \\
& = & \lambda ^{2}\left\Vert u_{3}\right\Vert ^{4}-\lambda ^{2}\left\Vert
u_{3}\right\Vert ^{4}=0.%
\end{array}%
\end{equation*}%
If we use $B(x_{1}+x_{2})+B(x_{1}-x_{2})=0,$ we get%
\begin{equation*}
f(x_{1},x_{1},x_{2},x_{2})+2f(x_{1},x_{2},x_{1},x_{2})=0.
\end{equation*}%
Changing $x_{2}$ into $x_{1}+x_{2}$, we have
\begin{equation*}
f(x_{1},x_{2},x_{2},x_{2})=0.
\end{equation*}%
In the last equation, if we change $x_{1}$ into $u_{2}$ and $x_{2}$ into $u_{1}$ where $%
u_{1}\bot u_{2},$ we obtain%
\begin{equation*}
f(u_{2},u_{1},u_{1},u_{1})=g_{_{N}}((\nabla \mathcal{T}_{\ast })(u_{2},u_{1}),(\nabla \mathcal{T}_{\ast })(u_{1},u_{1}))-\lambda
^{2}g_{_{M}}(u_{2},u_{1})g_{_{M}}(u_{1},u_{1})=0.
\end{equation*}%
So, we have%
\begin{equation*}
g_{_{\mathfrak{M}_2}}((\nabla \mathcal{T}_{\ast })(u_{2},u_{1}),(\nabla \mathcal{T}_{\ast })(u_{1},u_{1}))=0.
\end{equation*}
Conversely, we suppose that (3.1) is satisfied along the Riemannian map $\mathcal{T} .$
For an arbitrary orthogonal pair $x_{1},x_{2}\in \Gamma ((\ker \mathcal{T}_{\ast
p_{1}})^{\perp }),$ we have
\begin{equation*}
f(x_{1},x_{1},x_{1},x_{2})=g_{_{\mathfrak{M}_2}}((\nabla \mathcal{T}_{\ast })(x_{1},x_{1}),(\nabla \mathcal{T}_{\ast })(x_{1},x_{2}))-\lambda
^{2}g_{_{\mathfrak{M}_2}}(\mathcal{T}_{\ast }(x_{1}),\mathcal{T}_{\ast }(x_{1}))g_{_{\mathfrak{M}_2}}(\mathcal{T}_{\ast }(x_{1}),\mathcal{T}_{\ast }(x_{2}))=0.
\end{equation*}%
If we write $x_{1}+x_{2} $ instead of $x_{2}$, we have
\begin{equation*}
f(x_{1},x_{1},x_{1},x_{1})=0.
\end{equation*}
So we have
\begin{equation*}
\begin{array}{c}
g_{_{\mathfrak{M}_2}}((\nabla \mathcal{T}_{\ast })(x_{1},x_{1}),(\nabla \mathcal{T}_{\ast })(x_{1},x_{1}))=\lambda
^{2}g_{_{\mathfrak{M}_2}}(\mathcal{T}_{\ast }(x_{1}),\mathcal{T}_{\ast }(x_{1}))g_{_{N}}(\mathcal{T}_{\ast }(x_{1}),\mathcal{T}_{\ast }(x_{1})),%
\end{array}%
\end{equation*}%
that is, $\mathcal{T}$ is $\lambda h- $ isotropic.
\end{proof}
Let $(\mathfrak{M}_1,g_{_{\mathfrak{M}_1}})$ and $(\mathfrak{M}_2,g_{_{\mathfrak{M}_2}})$ be Riemannian manifolds and suppose
that $\mathcal{T}:(\mathfrak{M}_1,g_{_{\mathfrak{M}_1}})\rightarrow (\mathfrak{M}_2,g_{_{\mathfrak{M}_2}})$ is a Riemannian map between them.
Let $\alpha $ be a horizontal curve with curvature $\kappa $ in $\mathfrak{M}_1$ and $%
\gamma =\mathcal{T}\circ \alpha $ a curve with $\tilde{\kappa}$ in $\mathfrak{M}_2$ along $\mathcal{T}.$
Using (\ref{2.1}) and (\ref{2.2}), we obtain for $\forall t\in \mathbb{R}$%
\begin{equation*}
\begin{array}{ccl}
\tilde{\kappa}^{2} & = & g_{_{\mathfrak{M}_2}}(\overset{2}{\nabla _{\dot{\alpha}}^{\mathcal{T}}}%
\mathcal{T}_{\ast }(\dot{\alpha}),\overset{2}{\nabla _{\dot{\alpha}}^{\mathcal{T}}}\mathcal{T}_{\ast }(%
\dot{\alpha})), \\
& = & g_{_{\mathfrak{M}_2}}(\mathcal{T}_{\ast }(\overset{\mathfrak{M}_1}{\nabla }_{\dot{\alpha}}\dot{\alpha}%
)+(\nabla \mathcal{T}_{\ast })(\dot{\alpha},\dot{\alpha}),\mathcal{T}_{\ast }(\overset{\mathfrak{M}_1}{\nabla
}_{\dot{\alpha}}\dot{\alpha})+(\nabla \mathcal{T}_{\ast })(\dot{\alpha},\dot{\alpha})),
\\
& = & g_{_{\mathfrak{M}_1}}(\overset{\mathfrak{M}_1}{\nabla }_{\dot{\alpha}}\dot{\alpha},\overset{\mathfrak{M}_1}{%
\nabla }_{\dot{\alpha}}\dot{\alpha})+\left\Vert (\nabla \mathcal{T}_{\ast })(\dot{%
\alpha},\dot{\alpha})\right\Vert ^{2}, \\
& = & \kappa ^{2}+\left\Vert (\nabla \mathcal{T}_{\ast })(\dot{\alpha},\dot{\alpha}%
)\right\Vert ^{2},%
\end{array}%
\end{equation*}%
where $ \overset{\mathfrak{M}_2}{\nabla^\mathcal{T}} $ denotes the pullback connection of $ \overset{\mathfrak{M}_2}{\nabla} $. So, we can write%
\begin{equation}
\tilde{\kappa}=\tilde{\kappa}(t)=\sqrt{\kappa (t)^{2}+\left\Vert (\nabla \mathcal{T}_{\ast })(\dot{%
\alpha}(t),\dot{\alpha}(t))\right\Vert ^{2}}.  \label{3.1}
\end{equation}
We suppose that $ \mathcal{T}:\mathfrak{M}_1\rightarrow \mathfrak{M}_2 $ is a Riemannian map and $\alpha =\alpha (s)$ is a horizontal circle parametrized by arc length $s$ with
\begin{equation}
\overset{\mathfrak{M}_1}{\nabla} _{\dot{\alpha}}\dot{\alpha}=\kappa Y,\text{ \ }\overset{\mathfrak{M}_1}{\nabla} _{\dot{\alpha}}Y=-\kappa \dot{\alpha},
\label{circle}
\end{equation}%
where $Y=Y_{s}$ of unit vector along $ \alpha $ and $ \kappa $ is the curvature of $ \alpha $ on  $ \mathfrak{M}_1 $. For each point $ p\in \mathfrak{M}_1 $, each orthonormal couple $ u_{1},u_{2}\in\Gamma(ker\mathcal{T}_{*})^\bot $ at $ p $ and each constant $ \kappa>0 $, there is locally a unique horizontal circle $ \alpha=\alpha(s) $ on $ \mathfrak{M}_1 $ with initial condition that $ \alpha (0)=p, $ $ \dot{\alpha} (0)=u_{1} $ and  $ \overset{\mathfrak{M}_1}{\nabla}_{\dot{\alpha}} \dot{\alpha} (0)=\kappa u_{2} $ (see, \cite{Maeda, Sahin4}).
\begin{theorem}
Let $\mathcal{T}:\mathfrak{M}_1\rightarrow \mathfrak{M}_2$ be a smooth map between Riemannian manifolds $(\mathfrak{M}_1,g_{_{\mathfrak{M}_1}})$ and $(\mathfrak{M}_2,g_{_{\mathfrak{M}_2}})$. The following are equivalent:
\begin{enumerate}
  \item [(i)] $\mathcal{T}$ is a $ \lambda h-$isotropic Riemannian map,
  \item [(ii)] There is $ \kappa > 0 $ satisfying that for each horizontal circle $ \alpha $ with curvature $ \kappa $ on $ \mathfrak{M}_1 $, $ \gamma=\mathcal{T} \circ \alpha $ on $ \mathfrak{M}_2 $ has constant curvature $ \tilde{\kappa} $ along $ \gamma=\mathcal{T} \circ \alpha $.
\end{enumerate}
\end{theorem}
\begin{proof}
Assume that $ \mathcal{T}:\mathfrak{M}_1 \rightarrow \mathfrak{M}_2 $ is an isotropic Riemannin map. From (\ref{circle}) and (\ref{3.1}), curvature of $ \gamma $
\begin{equation}
\tilde{\kappa}(t)=\sqrt{\kappa^2+\lambda^2}
\end{equation}
is a constant.
Conversely, we suppose that there is $ \kappa > 0 $ satisfying that for each circle $ \alpha $ with curvature $ \kappa $ on $ \mathfrak{M}_1 $, the curve $ \gamma=\mathcal{T} \circ \alpha $ on $ \mathfrak{M}_2 $ has constant first curvature $ \tilde{\kappa} $ along this curve. Let $ u_{1},u_{2}\in\Gamma(ker\mathcal{T}_{*})^\bot $ be arbitrary orthonormal pair of vectors at $ p\in \mathfrak{M}_1 $. Suppose that $ \alpha=\alpha(s) $, $ s\in I $ be a circle with $ \kappa $ on $ \mathfrak{M}_1 $ with initial conditions $ \alpha (0)=p, $ $ \dot{\alpha} (0)=u_{1} $ and  $ \overset{\mathfrak{M}_1}{\nabla}_{\dot{\alpha}} \dot{\alpha} (0)=\kappa u_{2}. $ From (\ref{2.1}), we have
\begin{equation*}
 \tilde{\kappa }^{2}=\kappa ^{2}+\left\Vert (\nabla \mathcal{T}_{\ast })(\dot{\alpha},\dot{\alpha}%
)\right\Vert ^{2}.
\end{equation*}
Since $ \tilde{\kappa } $ is a constant and $ \alpha $ is a circle on $ \mathfrak{M}_1 $, $ \left \Vert (\nabla \mathcal{T}_{\ast })(\dot{\alpha},\dot{\alpha})\right\Vert $ is a constant on $ I $. Then,  we have
\begin{equation*}
0=\nabla _{\dot{\alpha}}^{\mathcal{T}}\left( g_{_{\mathfrak{M}_2}}((\nabla \mathcal{T}_{\ast })(\dot{\alpha},\dot{\alpha}),(\nabla
\mathcal{T}_{\ast })(\dot{\alpha},\dot{\alpha}))\right)
=2g_{_{\mathfrak{M}_2}}((\nabla _{\dot{\alpha}}^{\mathcal{T}})^{\bot}(\nabla \mathcal{T}_{\ast })(\dot{\alpha},\dot{\alpha}),(\nabla
\mathcal{T}_{\ast })(\dot{\alpha},\dot{\alpha})).
\end{equation*}
Using (\ref{circle}) and (\ref{2.6}), we get
\begin{equation}
g_{_{\mathfrak{M}_2}}(\tilde{\nabla}_{\dot{\alpha}}^{\mathcal{T}}(\nabla \mathcal{T}_{\ast })(\dot{\alpha},\dot{\alpha}),(\nabla
\mathcal{T}_{\ast })(\dot{\alpha},\dot{\alpha}))+2\kappa g_{_{\mathfrak{M}_2}}( (\nabla
\mathcal{T}_{\ast })(\dot{\alpha},Y_{s}),(\nabla
\mathcal{T}_{\ast })(\dot{\alpha},\dot{\alpha}))=0.
\label{3.8}
\end{equation}
Evaluating equation (\ref{3.8}) at $ s=0 $, we get
\begin{equation}
g_{_{\mathfrak{M}_2}}(\tilde{\nabla}_{u_{1}}^{\mathcal{T}}(\nabla \mathcal{T}_{\ast })(u_{1},u_{1}),(\nabla
\mathcal{T}_{\ast })(u_{1},u_{1}))+2\kappa g_{_{\mathfrak{M}_2}}( (\nabla
\mathcal{T}_{\ast })(u_{1},u_{2}),(\nabla
\mathcal{T}_{\ast })(u_{1},u_{2}))=0.
\label{3.9}
\end{equation}
Also, for another horizontal circle $ \beta=\beta(s) $ of the same curvature $ \kappa $ on $ \mathfrak{M}_1 $ with initial conditions $ \beta (0)=p, $ $ \dot{\beta} (0)=u_{1} $ and  $ \overset{1}{\nabla}_{\dot{\beta}} \dot{\beta} (0)=-\kappa u_{2}, $ we have
\begin{equation}
g_{_{\mathfrak{M}_2}}(\tilde{\nabla}_{u_{1}}^{\mathcal{T}}(\nabla \mathcal{T}_{\ast })(u_{1},u_{1}),(\nabla
\mathcal{T}_{\ast })(u_{1},u_{1}))-2\kappa g_{_{\mathfrak{M}_2}}( (\nabla
\mathcal{T}_{\ast })(u_{1},u_{2}),(\nabla
\mathcal{T}_{\ast })(u_{1},u_{1}))=0
\label{3.10}
\end{equation}
which corresponds to (\ref{3.9}). Then from (\ref{3.9}) and (\ref{3.10}), we obtain
\begin{equation*}
\kappa g_{_{\mathfrak{M}_2}}( (\nabla
\mathcal{T}_{\ast })(u_{1},u_{2}),(\nabla
\mathcal{T}_{\ast })(u_{1},u_{1}))=0.
\end{equation*}
Taking into consideration Lemma \ref{Lemma}, we can see that $ \mathcal{T} $ is a $ \lambda h-$ isotropic Riemannian map.
\end{proof}
\begin{corollary}
A totally umbilical Riemannian map $\mathcal{T}$ is a h-isotropic at the point $p_{1}$.
Conversely a h-isotropic Riemannian map $\mathcal{T}$ is a totally umbilical at $p_{1}$
if it satisfies%
\begin{equation*}
(\nabla \mathcal{T}_{\ast })(X,Y)=0
\end{equation*}%
for orthonormal vector fields $X$ and $Y\ $at $p_{1}$\ in $\Gamma ((\ker
\mathcal{T}_{\ast p_{1}})^{\perp }).$
\end{corollary}

\begin{proof}
Assume that $\mathcal{T}$ is a totally umbilical Riemannian map at the point $p_{1}$.
Then for $X_{1},X_{2}\in \Gamma ((\ker \mathcal{T}_{\ast p_{1}})^{\perp }),$ we have
\begin{equation*}
(\nabla \mathcal{T}_{\ast })(X,Y)=g_{_{\mathfrak{M}_1}}(X,Y)H_{2}.
\end{equation*}%
Especially
\begin{equation*}
(\nabla \mathcal{T}_{\ast })(X,Y)=0
\end{equation*}%
if $X$ and $Y$ are orthogonal. Then $\mathcal{T}$ is h-isotropic. Conversely we
suppose that $\mathcal{T}$ is a h-isotropic Riemannian map at $p_{1}$. Then  $\nabla \mathcal{T}_{\ast }$ satisfies $(3.4)$ for an arbitrary
orthogonal couple $X,Y\in \Gamma ((\ker \mathcal{T}_{\ast p_{1}})^{\perp }).$ Since $%
(\nabla \mathcal{T}_{\ast })(X,Y)=0$ and $\nabla \mathcal{T}_{\ast }$ is linear$,$ we have for
the orthogonal couple $\frac{1}{\sqrt{2}}(X+Y)$ and $\frac{1}{\sqrt{2}%
}(X-Y)$
\begin{equation*}
(\nabla \mathcal{T}_{\ast })(X,X)=(\nabla \mathcal{T}_{\ast })(Y,Y).
\end{equation*}%
Let $\left\{ X_{1},X_{2},...X_{n}\right\} $ denote an orthonormal frame in $%
\Gamma ((\ker \mathcal{T}_{\ast p_{1}})^{\perp }),$ then we get%
\begin{equation*}
(\nabla \mathcal{T}_{\ast })(X_{1},X_{1})=(\nabla \mathcal{T}_{\ast })(X_{2},X_{2})=...=(\nabla
\mathcal{T}_{\ast })(X_{n},X_{n}).
\end{equation*}%
Thus we have%
\begin{equation*}
H_{2}=(\nabla \mathcal{T}_{\ast })(X_{1},X_{1}),
\end{equation*}%
where $H_{2}$ is the mean curvature vector field of distribution $%
range\mathcal{T}_{\ast }.$ Moreover, choosing $X=\sum\limits_{i}a_{i}X_{i}$ and $%
Y=\sum\limits_{j}b_{j}X_{j},$ we get%
\begin{equation*}
0=(\nabla \mathcal{T}_{\ast })(X,Y)=\sum\limits_{i,j}a_{i}b_{j}(\nabla \mathcal{T}_{\ast
})(X_{i},X_{i})=g_{M}(X,Y)H_{2},
\end{equation*}%
which shows that $\mathcal{T}$ is umbilical.
\end{proof}

\section{A characterization of Riemannian maps in terms of helices}

We prove the following theorem which shows the effect of transforming helices to the base manifold  on the type of Riemannian maps in this section.
\begin{theorem}
Let $\mathcal{T}$ be a Riemannian map from a connected Riemannian manifold $%
(\mathfrak{M}_1,g_{_{\mathfrak{M}_1}}),$ $dim \mathfrak{M}_1 \geq 2$ to a Riemannian manifold $(\mathfrak{M}_2,g_{_{\mathfrak{M}_2}}).$ Let $%
\alpha $ be a horizontal helix with curvature $\kappa $ and torsion $\tau $
on $\mathfrak{M}_1,$ then $\mathcal{T}$ is umbilical and the mean curvature vector field $H_{2}$
satisfies the following equation%
\begin{equation*}
\left( \nabla _{\xi_{s}}^{\mathcal{T}^{\perp }}\right) ^{2}H_{2}=-\tau ^{2}H_{2},
\end{equation*}%
if and only if for every horizontal helix $\alpha $ on the base manifold $\mathfrak{M}_1$,  the corresponding  curve $\mathcal{T}\circ
\alpha $ is a helix on $\mathfrak{M}_2$.
\end{theorem}

\begin{proof}
We consider that $p\in \mathfrak{M}_1$\ and $\alpha (s)$ is a horizontal helix with curvature $\kappa $
and torsion $\tau $ on the base  manifold $\mathfrak{M}_1$, $\mathcal{T}\circ \alpha :I\rightarrow \mathfrak{M}_2$ is the corresponding
curve and we can define a
vector field $\mathcal{T}_{\ast }\xi$ along $\mathcal{T}\circ \alpha $ by%
\begin{equation*}
\mathcal{T}_{\ast }\xi(s))=\mathcal{T}_{\ast \alpha (s)}\xi(s),
\end{equation*}%
for each vector field $\xi_{s}$ along $\alpha $, where $\xi_{s}$ is the unit tangent vector field along $\alpha $ and $s$ is the arc length parameter. Now we assume that $\mathcal{T}\circ
\alpha $ is a helix with the curvature $\tilde{\kappa}$ and torsion $\tilde{%
\tau}$ on $\mathfrak{M}_2.$ From (\ref{helix}),  we have
\begin{equation}
\left( \overset{2~~~}{\nabla _{_{\xi_{s}}}^{\mathcal{T}}}\right) ^{3}\mathcal{T}_{\ast }(\xi_{s})+%
\tilde{K}^{2}\overset{2~~~~~}{\nabla _{_{\xi_{s}}}^{\mathcal{T}}}\mathcal{T}_{\ast }(\xi_{s})=0
\label{4.1}
\end{equation}%
where $\tilde{K}^{2}=\tilde{\kappa}^{2}+\tilde{\tau}^{2}.$ Using (\ref{2.1}), we
have%
\begin{equation}
\overset{2~~~}{\nabla _{_{\xi_{s}}}^{\mathcal{T}}}\mathcal{T}_{\ast }(\xi_{s})=\mathcal{T}_{\ast }({\overset{1}%
{\nabla }_{_{\xi_{s}}}\xi_{s}})+(\nabla \mathcal{T}_{\ast })(\xi_{s},\xi_{s}).
\label{4.2}
\end{equation}%
From (\ref{2.1}) and (\ref{4.2}), we obtain%
\begin{equation*}
\overset{2~~~}{\nabla _{_{\xi_{s}}}^{\mathcal{T}}}\mathcal{T}_{\ast }({\overset{1}{\nabla }%
_{_{\xi_{s}}}\xi_{s}})=\mathcal{T}_{\ast }((\overset{1}{\nabla }_{_{\xi_{s}}})^{2}\xi_{s})+(%
\nabla \mathcal{T}_{\ast })(\xi_{s},{\overset{1}{\nabla }_{_{\xi_{s}}}\xi_{s}}).
\end{equation*}%
Using (\ref{4.2}) and (\ref{2.3}) in the last equation, we get
\begin{eqnarray}
\left( \overset{2~~~}{\nabla _{_{\xi_{s}}}^{\mathcal{T}}}\right) ^{2}\mathcal{T}_{\ast }(\xi_{s})
&=&-S_{_{(\nabla \mathcal{T}_{\ast })(\xi_{s},\xi_{s})}}\mathcal{T}_{\ast }(\xi_{s})+\nabla
_{\xi_{s}}^{\mathcal{T}^{\perp }}(\nabla \mathcal{T}_{\ast })(\xi_{s},\xi_{s})+\mathcal{T}_{\ast }(({\overset{1}{\nabla }_{_{\xi_{s}}}})^{2}\xi_{s})  \notag \\
&&+(\nabla \mathcal{T}_{\ast })(\xi_{s},{\overset{1}{\nabla }_{_{\xi_{s}}}\xi_{s}}).
\label{4.3}
\end{eqnarray}%
Using (\ref{2.1}), (\ref{2.3}) and (\ref{4.3}), we get%
\begin{eqnarray*}
\left( \overset{2~~~}{\nabla _{_{\xi_{s}}}^{\mathcal{T}}}\right) ^{3}\mathcal{T}_{\ast }(\xi_{s})
&=&\mathcal{T}_{\ast }(\overset{1~~~}{\nabla _{_{\xi_{s}}}^{3}}\xi_{s})+(\nabla \mathcal{T}_{\ast
})(\xi_{s},((\overset{1}{\nabla }_{_{\xi_{s}}})^{2}\xi_{s}))-\overset{2~~~}{\nabla
_{_{\xi_{s}}}^{\mathcal{T}}}(S_{_{(\nabla \mathcal{T}_{\ast })(\xi_{s},\xi_{s})}}\mathcal{T}_{\ast }(\xi_{s})) \\
&&+\overset{2~~~~}{\nabla _{_{\xi_{s}}}^{\mathcal{T}}}(\nabla _{\xi_{s}}^{\mathcal{T}^{\perp
}}(\nabla \mathcal{T}_{\ast })(\xi_{s},\xi_{s}))+\overset{2~~~~}{\nabla _{_{\xi_{s}}}^{\mathcal{T}}}%
((\nabla \mathcal{T}_{\ast })(\xi_{s},{\overset{1}{\nabla }_{_{\xi_{s}}}\xi_{s}})).
\end{eqnarray*}%
So, we obtain%
\begin{eqnarray}
\left( \overset{2~~~}{\nabla _{_{\xi_{s}}}^{\mathcal{T}}}\right) ^{3}\mathcal{T}_{\ast }(\xi_{s})
&=&\mathcal{T}_{\ast }(\overset{1~~~}{\nabla _{_{\xi_{s}}}^{3}}\xi_{s})+(\nabla \mathcal{T}_{\ast
})(\xi_{s},(\overset{1}{\nabla }%
_{_{\xi_{s}}})^{2}\xi_{s})-\mathcal{T}_{\ast }(\overset{1}{%
\nabla }_{\xi_{s}}\text{~}^{\ast }\mathcal{T}_{\ast }S_{_{(\nabla \mathcal{T}_{\ast
})(\xi_{s},\xi_{s})}}\mathcal{T}_{\ast }(\xi_{s}))  \notag \\
&&-(\nabla \mathcal{T}_{\ast })(\xi_{s},S_{_{(\nabla \mathcal{T}_{\ast })(\xi_{s},\xi_{s})}}\mathcal{T}_{\ast
}(\xi_{s}))-S_{\nabla _{\xi_{s}}^{\mathcal{T}^{\perp }}(\nabla \mathcal{T}_{\ast
})(\xi_{s},\xi_{s})}\mathcal{T}_{\ast }(\xi_{s})  \notag \\
&&+\left( \nabla _{\xi_{s}}^{\mathcal{T}^{\perp }}\right) ^{2}(\nabla \mathcal{T}_{\ast
})(\xi_{s},\xi_{s})-S_{(\nabla \mathcal{T}_{\ast })(\xi_{s},{\overset{1}{\nabla }%
_{_{\xi_{s}}}\xi_{s}})}\mathcal{T}_{\ast }(\xi_{s}) \label{4.4} \\
&&+\nabla _{\xi_{s}}^{\mathcal{T}^{\perp }}(\nabla \mathcal{T}_{\ast })(\xi_{s},{\overset{1}{\nabla }%
_{_{\xi_{s}}}\xi_{s}}).  \notag
\end{eqnarray}%
Substituting (\ref{4.4}) and (\ref{4.2}) into (\ref{4.1}), we get%
\begin{equation}
\begin{array}{c}
\mathcal{T}_{\ast }(\overset{1~~~}{\nabla _{_{\xi_{s}}}^{3}}\xi_{s})+(\nabla \mathcal{T}_{\ast
})(\xi_{s},\overset{1~~~}{\nabla _{_{\xi_{s}}}^{2}}\xi_{s})-\mathcal{T}_{\ast }(\overset{1}{%
\nabla }_{\xi_{s}}\text{~}^{\ast }\mathcal{T}_{\ast }S_{_{(\nabla \mathcal{T}_{\ast
})(\xi_{s},\xi_{s})}}\mathcal{T}_{\ast }(\xi_{s})) \\
-(\nabla \mathcal{T}_{\ast })(\xi_{s},S_{_{(\nabla \mathcal{T}_{\ast })(\xi_{s},\xi_{s})}}\mathcal{T}_{\ast
}(\xi_{s}))-S_{\nabla _{\xi_{s}}^{\mathcal{T}^{\perp }}(\nabla \mathcal{T}_{\ast
})(\xi_{s},\xi_{s})}\mathcal{T}_{\ast }(\xi_{s}) \\
+\left( \nabla _{\xi_{s}}^{\mathcal{T}^{\perp }}\right) ^{2}(\nabla \mathcal{T}_{\ast
})(\xi_{s},\xi_{s})-S_{(\nabla \mathcal{T}_{\ast })(\xi_{s},{\overset{1}{\nabla }%
_{_{\xi_{s}}}\xi_{s}})}\mathcal{T}_{\ast }(\xi_{s}) \\
+\nabla _{\xi_{s}}^{\mathcal{T}^{\perp }}(\nabla \mathcal{T}_{\ast })(\xi_{s},{\overset{1}{\nabla }%
_{_{\xi_{s}}}\xi_{s}})+\tilde{K}^{2}\mathcal{T}_{\ast }({\overset{1}{\nabla }%
_{_{\xi_{s}}}\xi_{s}})+\tilde{K}^{2}(\nabla \mathcal{T}_{\ast })(\xi_{s},\xi_{s})=0.%
\end{array}%
\label{4.5}
\end{equation}%
By looking at $range\mathcal{T}_{\ast }$ and $(range\mathcal{T}_{\ast })^{\perp }$ components of (\ref{4.5}), we have%
\begin{equation*}
\begin{array}{c}
\mathcal{T}_{\ast }(\overset{1~~}{\nabla _{_{\xi_{s}}}^{3}}\xi_{s})+\tilde{K}^{2}\mathcal{T}_{\ast }(%
{\overset{1}{\nabla }_{_{\xi_{s}}}\xi_{s}})-\mathcal{T}_{\ast }(\overset{1}{\nabla }%
_{\xi_{s}}\text{~}^{\ast }\mathcal{T}_{\ast }S_{_{(\nabla \mathcal{T}_{\ast
})(\xi_{s},\xi_{s})}}\mathcal{T}_{\ast }(\xi_{s})) \\
-S_{\nabla _{\xi_{s}}^{\mathcal{T}^{\perp }}(\nabla \mathcal{T}_{\ast })(\xi_{s},\xi_{s})}\mathcal{T}_{\ast
}(\xi_{s})-S_{(\nabla \mathcal{T}_{\ast })(\xi_{s},{\overset{1}{\nabla }_{_{\xi_{s}}}\xi_{s}}%
)}\mathcal{T}_{\ast }(\xi_{s})=0%
\end{array}%
\end{equation*}%
and%
\begin{equation*}
\begin{array}{c}
(\nabla \mathcal{T}_{\ast })(\xi_{s},(\overset{1}{\nabla }_{_{\xi_{s}}})^{2}\xi_{s})-(\nabla
\mathcal{T}_{\ast })(\xi_{s},S_{_{(\nabla \mathcal{T}_{\ast })(\xi_{s},\xi_{s})}}\mathcal{T}_{\ast
}(\xi_{s}))+\left( \nabla _{\xi_{s}}^{\mathcal{T}^{\perp }}\right) ^{2}(\nabla \mathcal{T}_{\ast
})(\xi_{s},\xi_{s}) \\
+\nabla _{\xi_{s}^{\mathcal{T}^{\perp }}(\nabla \mathcal{T}_{\ast })(\xi_{s},{\overset{1}{\nabla }}%
_{_{\xi_{s}}}\xi_{s}})+\tilde{K}^{2}(\nabla \mathcal{T}_{\ast })(\xi_{s},\xi_{s})=0.%
\end{array}%
\end{equation*}%
Considering (\ref{2.6}),we arrive at
\begin{eqnarray*}
\left( \tilde{\nabla}_{\xi_{s}}S\right) _{(\nabla \mathcal{T}_{\ast })(\xi_{s},\xi_{s})}\mathcal{T}_{\ast
}(\xi_{s}) &=&\mathcal{T}_{\ast }(\overset{1}{\nabla }_{\xi_{s}}\text{~}^{\ast }\mathcal{T}_{\ast
}S_{_{(\nabla \mathcal{T}_{\ast })(\xi_{s},\xi_{s})}}\mathcal{T}_{\ast }(\xi_{s})) \\
&&-S_{_{\nabla _{\xi_{s}}^{\mathcal{T}^{\perp }}(\nabla \mathcal{T}_{\ast })(\xi_{s},\xi_{s})}}\mathcal{T}_{\ast
}(\xi_{s})-S_{(\nabla \mathcal{T}_{\ast })(\xi_{s},\xi_{s})}\overset{2}{P\nabla _{\xi_{s}}^{\mathcal{T}}}\mathcal{T}_{\ast
}(\xi_{s})
\end{eqnarray*}%
and%
\begin{eqnarray}
\mathcal{T}_{\ast }(\overset{1}{\nabla }_{\xi_{s}}\text{~}^{\ast }\mathcal{T}_{\ast }S_{_{(\nabla
\mathcal{T}_{\ast })(\xi_{s},\xi_{s})}}\mathcal{T}_{\ast }(\xi_{s})) &=&\left( \tilde{\nabla}%
_{\xi}S\right) _{(\nabla \mathcal{T}_{\ast })(\xi_{s},\xi_{s})}\mathcal{T}_{\ast }(\xi_{s})  \notag \\
&&+S_{_{\nabla _{\xi_{s}}^{\mathcal{T}^{\perp }}(\nabla \mathcal{T}_{\ast })(\xi_{s},\xi_{s})}}\mathcal{T}_{\ast
}(\xi_{s})+S_{(\nabla \mathcal{T}_{\ast })(\xi_{s},\xi_{s})}\overset{2}{P\nabla _{\xi_{s}}^{\mathcal{T}}}\mathcal{T}_{\ast
}(\xi_{s}).
\label{4.6}
\end{eqnarray}%
Using (\ref{4.6}) and Frenet formulas in the tangent part of (\ref{4.5}) we get
\begin{equation*}
\begin{array}{c}
\mathcal{T}_{\ast }(-\kappa (\kappa ^{2}+\tau ^{2})V_2)+\tilde{K}^{2}\mathcal{T}_{\ast }(\kappa
V_2)-\left( \tilde{\nabla}_{\xi_{s}}S\right) _{(\nabla \mathcal{T}_{\ast
})(\xi_{s},\xi_{s})}\mathcal{T}_{\ast }(\xi_{s})-2S_{_{\nabla _{\xi_{s}}^{\mathcal{T}^{\perp }}(\nabla
\mathcal{T}_{\ast })(\xi_{s},\xi_{s})}}\mathcal{T}_{\ast }(\xi_{s}) \\
-S_{(\nabla \mathcal{T}_{\ast })(\xi_{s},\xi_{s})}\overset{2}{P\nabla _{\xi_{s}}^{\mathcal{T}}}\mathcal{T}_{\ast
}(\xi_{s})-KS_{_{(\nabla \mathcal{T}_{\ast })(\xi_{s},V_2)}}\mathcal{T}_{\ast }(\xi_{s})=0.%
\end{array}%
\end{equation*}%
From the last equation, we can write %
\begin{equation}
\begin{array}{c}
(\tilde{K}^{2}-\kappa ^{2}-\tau ^{2})\kappa \mathcal{T}_{\ast }(V_2)-\kappa S_{_{(\nabla
\mathcal{T}_{\ast })(\xi_{s},V_2)}}\mathcal{T}_{\ast }(\xi_{s})-\left( \tilde{\nabla}_{\xi_{s}}S\right)
_{(\nabla \mathcal{T}_{\ast })(\xi_{s},\xi_{s})}\mathcal{T}_{\ast }(\xi_{s}) \\
-2S_{_{\nabla _{\xi_{s}}^{\mathcal{T}^{\perp }}(\nabla \mathcal{T}_{\ast })(\xi_{s},\xi_{s})}}\mathcal{T}_{\ast
}(\xi_{s})-S_{(\nabla \mathcal{T}_{\ast })(\xi_{s},\xi_{s})}\overset{2}{P\nabla _{\xi_{s}}^{\mathcal{T}}}\mathcal{T}_{\ast
}(\xi_{s})=0.%
\label{4.7}
\end{array}%
\end{equation}%
From (\ref{2.5}), we have%
\begin{equation*}
\left( \tilde{\nabla}_{\xi_{s}}(\nabla \mathcal{T}_{\ast })\right) (\xi_{s},\xi_{s})=\nabla
_{\xi_{s}}^{\mathcal{T}^{\perp }}(\nabla \mathcal{T}_{\ast })(\xi_{s},\xi_{s})-2(\nabla \mathcal{T}_{\ast
})(\xi_{s},{\overset{1}{\nabla }_{_{\xi_{s}}}\xi_{s}})
\end{equation*}%
or%
\begin{equation}
2\nabla _{\xi_{s}}^{\mathcal{T}^{\perp }}(\nabla \mathcal{T}_{\ast })(\xi_{s},\xi_{s})=2\left( \tilde{%
\nabla}_{\xi_{s}}(\nabla \mathcal{T}_{\ast })\right) (\xi_{s},\xi_{s})+4(\nabla \mathcal{T}_{\ast
})(\xi_{s},{\overset{1}{\nabla }_{_{\xi_{s}}}\xi_{s}}).
\label{4.8}
\end{equation}%
Substituting (\ref{4.8}) into (\ref{4.7}), we have%
\begin{equation*}
\begin{array}{l}
(\tilde{K}^{2}-\kappa ^{2}-\tau ^{2})\kappa \mathcal{T}_{\ast }(V_2)-\kappa S_{_{(\nabla
\mathcal{T}_{\ast })(\xi_{s},V_2)}}\mathcal{T}_{\ast }(\xi_{s})-\left( \tilde{\nabla}_{\xi_{s}}S\right)
_{(\nabla \mathcal{T}_{\ast })(\xi_{s},\xi_{s})}\mathcal{T}_{\ast }(\xi_{s}) \\
-2S_{_{(\tilde{\nabla}_{\xi_{s}}(\nabla \mathcal{T}_{\ast }))(\xi_{s},\xi_{s})}}\mathcal{T}_{\ast
}(\xi_{s})-4\kappa S_{(\nabla \mathcal{T}_{\ast })(\xi_{s},V_2)}\mathcal{T}_{\ast }(\xi_{s})-S_{(\nabla \mathcal{T}_{\ast })(\xi_{s},\xi_{s})}\overset{2}{P\nabla _{\xi_{s}}^{\mathcal{T}}}\mathcal{T}_{\ast }(\xi_{s})=0.%
\end{array}%
\end{equation*}%
or%
\begin{eqnarray*}
&&(\tilde{K}^{2}-\kappa ^{2}-\tau ^{2})\kappa \mathcal{T}_{\ast }(V_2)-5\kappa
S_{(\nabla \mathcal{T}_{\ast })(\xi_{s},V_2)}\mathcal{T}_{\ast }(\xi_{s})-S_{(\nabla \mathcal{T}_{\ast })(\xi_{s},\xi_{s})}\overset{2}{P\nabla _{\xi_{s}}^{\mathcal{T}}}\mathcal{T}_{\ast }(\xi_{s}) \\
&=&2S_{_{(\tilde{\nabla}_{\xi_{s}}(\nabla \mathcal{T}_{\ast }))(\xi_{s},\xi_{s})}}\mathcal{T}_{\ast
}(\xi_{s})+\left( \tilde{\nabla}_{\xi_{s}}S\right) _{(\nabla \mathcal{T}_{\ast
})(\xi_{s},\xi_{s})}\mathcal{T}_{\ast }(\xi_{s}).
\end{eqnarray*}%
Taking inner product with $\mathcal{T}_{\ast }(\xi_{s}),$ we obtain%
\begin{equation}
\begin{array}{l}
(\tilde{K}^{2}-\kappa ^{2}-\tau ^{2})\kappa g_{_{\mathfrak{M}_2}}(\mathcal{T}_{\ast }(V_2),\mathcal{T}_{\ast
}(\xi_{s}))-5\kappa g_{_{\mathfrak{M}_2}}(S_{_{(\nabla \mathcal{T}_{\ast })(\xi_{s},V_2)}}\mathcal{T}_{\ast
}(\xi_{s}),\mathcal{T}_{\ast }(\xi_{s})) \\
-g_{_{\mathfrak{M}_2}}(S_{(\nabla \mathcal{T}_{\ast })(\xi_{s},\xi_{s})}\overset{2}{P\nabla _{\xi_{s}}^{\mathcal{T}}}\mathcal{T}_{\ast }(\xi_{s}),\mathcal{T}_{\ast }(\xi_{s}))=g_{_{\mathfrak{M}_2}}(2S_{_{(\tilde{\nabla}_{\xi_{s}}(\nabla
\mathcal{T}_{\ast }))(\xi_{s},\xi_{s})}}\mathcal{T}_{\ast }(X_{s}),\mathcal{T}_{\ast }(\xi_{s})) \\
\text{ \ \ \ \ \ \ \ \ \ \ \ \ \ \ \ \ \ \ \ \ \ \ \ \ \ \ \ \ \ \ \ \ \ \ \
\ \ \ \ \ \ \ \ \ \ \ \ \ \ \ \ \ \ \ \ \ }+g_{_{\mathfrak{M}_2}}(\left( \tilde{\nabla}%
_{\xi_{s}}S\right) _{(\nabla \mathcal{T}_{\ast })(\xi_{s},\xi_{s})}\mathcal{T}_{\ast }(\xi_{s}),\mathcal{T}_{\ast
}(\xi_{s})).%
\label{4.9}
\end{array}%
\end{equation}%
Since $\alpha $ is a horizontal curve, we have%
\begin{equation}
g_{_{\mathfrak{M}_2}}(\mathcal{T}_{\ast }(\xi_{s}),\mathcal{T}_{\ast }({\overset{1}{\nabla }_{_{\xi_{s}}}\xi_{s}}%
))=g_{_{\mathfrak{M}_1}}(\xi_{s},{\overset{1}{\nabla }_{_{\xi_{s}}}\xi_{s}})=\kappa
g_{_{\mathfrak{M}_1}}(\xi_{s},\kappa V_2)=0.
\label{4.10}
\end{equation}%
Taking into consideration (\ref{4.10}) and (\ref{2.4}), (\ref{4.9}) reduces to%
\begin{equation*}
\begin{array}{l}
-5\kappa g_{_{\mathfrak{M}_2}}((\nabla \mathcal{T}_{\ast })(\xi_{s},V_2),(\nabla \mathcal{T}_{\ast
})(\xi_{s},\xi_{s}))-g_{_{\mathfrak{M}_2}}(S_{(\nabla \mathcal{T}_{\ast })(\xi_{s},\xi_{s})}\overset{2}{P\nabla _{\xi_{s}}^{\mathcal{T}}}\mathcal{T}_{\ast }(\xi_{s}),\mathcal{T}_{\ast }(\xi_{s})) \\
=2g_{_{\mathfrak{M}_2}}((\tilde{\nabla}_{\xi_{s}}(\nabla \mathcal{T}_{\ast }))(\xi_{s},\xi_{s}),(\nabla
\mathcal{T}_{\ast })(\xi_{s},\xi_{s}))+g_{_{\mathfrak{M}_2}}(\left( \tilde{\nabla}_{\xi_{s}}S\right)
_{(\nabla \mathcal{T}_{\ast })(\xi_{s},\xi_{s})}\mathcal{T}_{\ast }(\xi_{s}),\mathcal{T}_{\ast }(\xi_{s})).%
\end{array}%
\end{equation*}%
Using (\ref{2.7}), we get%
\begin{equation*}
\begin{array}{l}
-5\kappa g_{_{\mathfrak{M}_2}}((\nabla \mathcal{T}_{\ast })(\xi_{s},V_2),(\nabla \mathcal{T}_{\ast
})(\xi_{s},\xi_{s}))-g_{_{\mathfrak{M}_2}}((\nabla \mathcal{T}_{\ast })(\xi_{s},{\overset{1}{\nabla }%
_{_{\xi_{s}}}\xi_{s}}),(\nabla \mathcal{T}_{\ast })(\xi_{s},\xi_{s})) \\
=3g_{_{\mathfrak{M}_2}}((\tilde{\nabla}_{\xi_{s}}(\nabla \mathcal{T}_{\ast }))(\xi_{s},\xi_{s}),(\nabla
\mathcal{T}_{\ast })(\xi_{s},\xi_{s})).%
\end{array}%
\end{equation*}%
By using (\ref{4.8}), we obtain%
\begin{equation*}
\begin{array}{l}
-6\kappa g_{_{\mathfrak{M}_2}}((\nabla \mathcal{T}_{\ast })(\xi_{s},V_2),(\nabla \mathcal{T}_{\ast
})(\xi_{s},\xi_{s}))=3g_{_{\mathfrak{M}_2}}(\nabla _{\xi_{s}}^{\mathcal{T}^{\perp }}(\nabla \mathcal{T}_{\ast
})(\xi_{s},\xi_{s}),(\nabla \mathcal{T}_{\ast })(\xi_{s},\xi_{s})) \\
\text{ \ \ \ \ \ \ \ \ \ \ \ \ \ \ \ \ \ \ \ \ \ \ \ \ \ \ \ \ \ \ \ \ \ \ \
\ \ \ \ \ \ \ \ \ \ \ \ \ \ \ \ \ \ \ \ \ \ \ \ \ \ \ }-6\kappa
g_{_{\mathfrak{M}_2}}((\nabla \mathcal{T}_{\ast })(\xi_{s},V_2),(\nabla \mathcal{T}_{\ast })(\xi_{s},\xi_{s})).%
\end{array}%
\end{equation*}%
Hence, we get%
\begin{equation*}
g_{_{\mathfrak{M}_2}}(\nabla _{\xi_{s}}^{\mathcal{T}^{\perp }}(\nabla \mathcal{T}_{\ast })(\xi_{s},\xi_{s}),(\nabla
\mathcal{T}_{\ast })(\xi_{s},\xi_{s}))=0
\end{equation*}%
or%
\begin{equation*}
\nabla _{\xi_{s}}^{\mathcal{T}}\left( g_{_{\mathfrak{M}_2}}((\nabla \mathcal{T}_{\ast })(\xi_{s},\xi_{s}),(\nabla
\mathcal{T}_{\ast })(\xi_{s},\xi_{s}))\right) =0.
\end{equation*}%
So, we can see that%
\begin{equation*}
\left\Vert (\nabla \mathcal{T}_{\ast })(\xi_{s},\xi_{s})\right\Vert =const..
\end{equation*}%
This implies that $\mathcal{T}$ is an isotropic Riemannian map. So, from Lemma $3.1$,
we have%
\begin{equation}
g_{_{\mathfrak{M}_2}}((\nabla \mathcal{T}_{\ast })(\xi_{s},\xi_{s}),(\nabla \mathcal{T}_{\ast })(\xi_{s},V_2))=0.
\label{4.11}
\end{equation}%
The normal part of (\ref{4.5}) reduces to%
\begin{equation}
(\tilde{K}^{2}-\kappa ^{2}-\tau ^{2})\kappa \mathcal{T}_{\ast }(V_2)=5\kappa
S_{_{(\nabla \mathcal{T}_{\ast })(\xi_{s},V_2)}}\mathcal{T}_{\ast }(\xi_{s})+S_{(\nabla \mathcal{T}_{\ast })(\xi_{s},\xi_{s})}\overset{2}{P\nabla _{\xi_{s}}^{\mathcal{T}}}\mathcal{T}_{\ast }(\xi_{s}).
\label{4.12}
\end{equation}%
Taking  derivative of (\ref{4.8}) and using the normal part of (\ref{4.5}) and (\ref{2.6}), the normal part of (\ref{4.5}) reduces to%
\begin{equation}
\begin{array}{c}
5\kappa (\tilde{\nabla}_{\xi_{s}}(\nabla \mathcal{T}_{\ast }))(\xi_{s},V_2)+3\kappa
^{2}(\nabla \mathcal{T}_{\ast })(V_2,V_2)+4\kappa \tau (\nabla \mathcal{T}_{\ast })(\xi_{s},V_3) \\
=(\nabla \mathcal{T}_{\ast })(\xi_{s},S_{(\nabla \mathcal{T}_{\ast })(\xi_{s},\xi_{s})}\mathcal{T}_{\ast
}(\xi_{s}))+4\kappa ^{2}(\nabla \mathcal{T}_{\ast })(\xi_{s},\xi_{s}) \\
-\tilde{\nabla}_{\xi_{s}}^{2}(\nabla \mathcal{T}_{\ast }))(\xi_{s},\xi_{s})-\tilde{K}%
^{2}(\nabla \mathcal{T}_{\ast })(\xi_{s},\xi_{s}).%
\label{4.13}
\end{array}%
\end{equation}%
Changing $V_3$ into $-V_3$ into (\ref{4.13}) we obtain%
\begin{equation}
(\nabla \mathcal{T}_{\ast })(\xi_{s},V_3)=0.
\label{4.14}
\end{equation}%
From (\ref{4.11}), (\ref{4.14}) and Corollary 3.1, we have that $\mathcal{T}$ is umbilical
map. If we change $V_2$ with $-V_2$ into (\ref{4.13}), we get%
\begin{equation}
5(\tilde{\nabla}_{\xi_{s}}(\nabla \mathcal{T}_{\ast }))(\xi_{s},V_2)=0.
\label{4.15}
\end{equation}%
Taking inner product (\ref{4.12}) with $\mathcal{T}_{\ast }(V_2),$ we have
\begin{eqnarray*}
(\tilde{K}^{2}-\kappa ^{2}-\tau ^{2})\kappa &=&5\kappa g_{_{\mathfrak{M}_2}}((\nabla
\mathcal{T}_{\ast })(\xi_{s},V_2),(\nabla \mathcal{T}_{\ast })(\xi_{s},V_2)) \\
&&+\kappa g_{_{\mathfrak{M}_2}}((\nabla \mathcal{T}_{\ast })(\xi_{s},\xi_{s}),(\nabla \mathcal{T}_{\ast })(V_2,V_2)),
\end{eqnarray*}%
that is%
\begin{equation}
(\tilde{K}^{2}-\kappa ^{2}-\tau ^{2})=\left\Vert H_{2}\right\Vert ^{2}.
\label{4.16}
\end{equation}%
Substituting (\ref{4.14}) and (\ref{4.15}) into (\ref{4.13}), we get%
\begin{equation*}
\begin{array}{c}
3\kappa ^{2}(\nabla \mathcal{T}_{\ast })(V_2,V_2)=(\nabla \mathcal{T}_{\ast }))(\xi_{s},S_{(\nabla
\mathcal{T}_{\ast })(\xi_{s},\xi_{s})}\mathcal{T}_{\ast }(\xi_{s})) \\
+4\kappa ^{2}(\nabla \mathcal{T}_{\ast })(\xi_{s},\xi_{s})-\tilde{\nabla}%
_{\xi_{s}}^{2}(\nabla \mathcal{T}_{\ast }))(\xi_{s},\xi_{s})-\tilde{K}^{2}(\nabla \mathcal{T}_{\ast
})(\xi_{s},\xi_{s}).%
\end{array}%
\end{equation*}%
From umbilicity and derivative of (\ref{4.8}), we have%
\begin{equation*}
\begin{array}{c}
3\kappa ^{2}H_{2}=(\nabla \mathcal{T}_{\ast }))(\xi_{s},S_{(\nabla \mathcal{T}_{\ast
})(\xi_{s},\xi_{s})}\mathcal{T}_{\ast }(\xi_{s}))+4\kappa ^{2}H_{2}-(\nabla
_{\xi_{s}}^{\mathcal{T}^{\perp }})^{2}H_{2} \\
+2\kappa ^{2}H_{2}-2\kappa ^{2}H_{2}-\tilde{K}^{2}H_{2}.%
\end{array}%
\end{equation*}%
From (\ref{2.4}), we can see that%
\begin{equation*}
(\tilde{K}^{2}-\kappa ^{2})H_{2}=g_{_{\mathfrak{M}_2}}((\nabla \mathcal{T}_{\ast
})(\xi_{s},\xi_{s}),(\nabla \mathcal{T}_{\ast })(\xi_{s},\xi_{s}))H_{2}-(\nabla
_{\xi_{s}}^{\mathcal{T}^{\perp }})^{2}H_{2}.
\end{equation*}%
So, we get%
\begin{equation*}
(\tilde{K}^{2}-\kappa ^{2}-\left\Vert H_{2}\right\Vert ^{2})H_{2}=-(\nabla
_{\xi_{s}}^{\mathcal{T}^{\perp }})^{2}H_{2}.
\end{equation*}%
Using (\ref{4.16}), we have%
\begin{equation*}
(\nabla _{\xi_{s}}^{\mathcal{T}^{\perp }})^{2}H_{2}=-\tau ^{2}H_{2}.
\end{equation*}%
Conversely, we assume that $\mathcal{T}$ is a umbilical map and mean curvature vector
field satisfies $\left( \nabla _{\xi_{s}}^{\mathcal{T}^{\perp }}\right) ^{2}H_{2}=-\tau
^{2}H_{2}.$ Then we calculate%
\begin{eqnarray*}
\left( \overset{2~~}{\nabla _{_{\xi_{s}}}^{\mathcal{T}}}\right) ^{3}\mathcal{T}_{\ast }(\xi_{s})
&=&\mathcal{T}_{\ast }(({\overset{1}{\nabla }{_{\xi_{s}}})^{3}\xi_{s}%
})+g_{_{\mathfrak{M}_2}}(\xi_{s},{(\overset{1}{\nabla }_{_{\xi_{s}}})^{2}\xi_{s}})H_{2}-\mathcal{T}_{\ast }(%
\overset{1}{\nabla }_{\xi_{s}}\text{~}^{\ast }\mathcal{T}_{\ast }(S_{_{(\nabla \mathcal{T}_{\ast
})(\xi_{s},\xi_{s})}}\mathcal{T}_{\ast }(\xi_{s}))) \\
&&-g_{_{\mathfrak{M}_2}}(\xi_{s},S_{_{g_{_{_{\mathfrak{M}_2}}}(\xi_{s},\xi_{s})H_{2}}}\mathcal{T}_{\ast
}(\xi_{s}))H_{2}-S_{_{\nabla _{\xi_{s}}^{\mathcal{T}^{\perp
}}g_{_{_{\mathfrak{M}_2}}}(\xi_{s},\xi_{s})H_{2}}}\mathcal{T}_{\ast }(\xi_{s}) \\
&&+\left( \nabla _{\xi_{s}}^{\mathcal{T}^{\perp }}\right)
^{2}g_{_{_{\mathfrak{M}_2}}}(\xi_{s},\xi_{s})H_{2}-S_{g_{_{\mathfrak{M}_2}}(\xi_{s},({{\overset{1}\nabla}_{{\xi_{s}}} )^{2}\xi_{s}})H_{2}}\mathcal{T}_{\ast }(\xi_{s})
\end{eqnarray*}%
or%
\begin{eqnarray}
\left( \overset{2~~}{\nabla _{_{\xi_{s}}}^{\mathcal{T}}}\right) ^{3}\mathcal{T}_{\ast }(\xi_{s})
&=&\mathcal{T}_{\ast }(({\overset{1}{\nabla}{_{\xi_{s}}})^{3} \xi_{s}})-\mathcal{T}_{\ast }(\overset{1%
}{\nabla }_{\xi_{s}}\text{~}^{\ast }\mathcal{T}_{\ast }(S_{_{H_{2}}}\mathcal{T}_{\ast }(\xi_{s})))  \label{4.17} \\
&&-(K^{2}+\left\Vert H_{2}\right\Vert ^{2})H_{2}-S_{_{\nabla
_{\xi_{s}}^{\mathcal{T}^{\perp }}H_{2}}}\mathcal{T}_{\ast }(\xi_{s})+\left( \nabla
_{\xi_{s}}^{\mathcal{T}^{\perp }}\right) ^{2}H_{2}.  \notag
\end{eqnarray}%
If we use (\ref{4.17}), we have%
\begin{equation*}
\begin{array}{c}
\left( \overset{2~~}{\nabla _{_{\xi_{s}}}^{\mathcal{T}}}\right) ^{3}\mathcal{T}_{\ast
}(\xi_{s})+(\kappa ^{2}+\tau ^{2}+\left\Vert H_{2}\right\Vert ^{2})\overset{2~~%
}{\nabla _{_{\xi_{s}}}^{\mathcal{T}}}\mathcal{T}_{\ast }(\xi_{s})=\mathcal{T}_{\ast }(({\overset{1}{\nabla}{_{\xi_{s}}})^{3} \xi_{s}}) \\
-\mathcal{T}_{\ast }(\overset{1}{\nabla }_{\xi_{s}}\text{~}^{\ast }\mathcal{T}_{\ast
}(S_{_{H_{2}}}\mathcal{T}_{\ast }(\xi_{s})))-(K^{2}+\left\Vert H_{2}\right\Vert
^{2})H_{2}-S_{_{\nabla _{\xi_{s}}^{\mathcal{T}^{\perp }}H_{2}}}\mathcal{T}_{\ast }(\xi_{s})+\left(
\nabla _{\xi_{S}}^{\mathcal{T}^{\perp }}\right) ^{2}H_{2} \\
+(\kappa ^{2}+\tau ^{2}+\left\Vert H_{2}\right\Vert ^{2})\mathcal{T}_{\ast }(\overset{%
1~~~~~~~~}{\nabla _{_{\xi_{s}}}\xi_{s}})+(\kappa ^{2}+\tau ^{2}+\left\Vert
H_{2}\right\Vert ^{2})H_{2}.%
\end{array}%
\end{equation*}%
or%
\begin{equation*}
\begin{array}{c}
\left( \overset{2~~}{\nabla _{_{\xi_{s}}}^{\mathcal{T}}}\right) ^{3}\mathcal{T}_{\ast
}(\xi_{s})+(\kappa ^{2}+\tau ^{2}+\left\Vert H_{2}\right\Vert ^{2})\overset{2~~%
}{\nabla _{_{\xi_{s}}}^{\mathcal{T}}}\mathcal{T}_{\ast }(\xi_{s})=\mathcal{T}_{\ast }(({\overset{1}{\nabla}{_{\xi_{s}}})^{3} \xi_{s}}) \\
-\mathcal{T}_{\ast }(\overset{1}{\nabla }_{\xi_{s}}\text{~}^{\ast }\mathcal{T}_{\ast
}(S_{_{H_{2}}}\mathcal{T}_{\ast }(\xi_{s})))-S_{_{\nabla _{\xi_{s}}^{\mathcal{T}^{\perp
}}H_{2}}}\mathcal{T}_{\ast }(\xi_{s})+(\kappa ^{2}+\tau ^{2}+\left\Vert H_{2}\right\Vert
^{2})\mathcal{T}_{\ast }({\overset{1}{\nabla }_{_{\xi_{s}}}\xi_{s}}).%
\end{array}%
\end{equation*}%
Then we get%
\begin{equation*}
\mathcal{T}_{\ast }(\left( {\overset{1}{\nabla}_{_{\xi_{s}} }}\right) ^{3}{\xi_{s}}%
)+(\kappa ^{2}+\tau ^{2})\mathcal{T}_{\ast }({\overset{1}{\nabla }_{_{\xi_{s}}}\xi_{s}})=0.
\end{equation*}%
Since $\alpha $ is a horizontal helix on $\mathfrak{M}_1$, then we find $\gamma =\mathcal{T}\circ
\alpha $ is a helix on $\mathfrak{M}_2$.
\end{proof}

\noindent{\bf Acknowledgement:}This work is supported by TUBITAK (The Scientific and Technological Council of Turkey) with project number 119F025.

\end{document}